\documentclass[graybox]{svmult}

\usepackage{graphics,epsf,epsfig,psfrag}

\usepackage{mathptmx}       
\usepackage{helvet}         
\usepackage{courier}        
\usepackage{type1cm}        
\usepackage{makeidx}         
\usepackage{graphicx}        
\usepackage{multicol}        
\usepackage[bottom]{footmisc}
\usepackage{color}

\makeindex             

\usepackage{amsmath, amsfonts}

\usepackage{hyperref}

\DeclareMathOperator{\sign}{\mathrm{sign}}

\newcommand{\ind}{\mathbf{1}}

\newcommand{\R}{\mathbb{R}}
\newcommand{\Z}{\mathbb{Z}}
\newcommand{\N}{\mathbb{N}}
\renewcommand{\tilde}{\widetilde}

\newcommand{\cN}{{\ensuremath{\mathcal N}} }
\newcommand{\cL}{{\ensuremath{\mathcal L}} }

\newcommand{\cD}{{\ensuremath{\mathcal D}} }

\newcommand{\bP}{{\ensuremath{\mathbf P}} }
\newcommand{\bE}{{\ensuremath{\mathbf E}} }

\DeclareMathSymbol{\leqslant}{\mathalpha}{AMSa}{"36} 
\DeclareMathSymbol{\geqslant}{\mathalpha}{AMSa}{"3E} 
\DeclareMathSymbol{\eset}{\mathalpha}{AMSb}{"3F}     
\newcommand{\dd}{\,\text{\rm d}}             

\newcommand{\bbE}{{\ensuremath{\mathbb E}} }

\newcommand{\bbP}{{\ensuremath{\mathbb P}} }

\newcommand{\bbZ}{{\ensuremath{\mathbb Z}} }

\newcommand{\ga}{\alpha}
\newcommand{\gb}{\beta}
\newcommand{\gd}{\delta}
\newcommand{\gep}{\varepsilon}       

\newcommand{\gD}{\Delta}

\newcommand{\go}{\omega}

\newcommand{\gO}{\Omega}
\newcommand{\gl}{\lambda}

\newcommand{\tf}{\textsc{f}}
\newcommand{\M}{\textsc{M}}

\newcommand{\rc}{\mathtt c}

\newcommand{\C}{\textsc{c}}

\begin{document}

\title*{Copolymers at selective interfaces: settled issues and open problems}
\author{Francesco Caravenna, Giambattista Giacomin and Fabio Lucio Toninelli}
\institute{Francesco Caravenna \at Dipartimento di Matematica Pura e Applicata, 
Universit\`a degli Studi di Padova,
via Trie\-ste 63, 35121 Padova, Italy,
\email{francesco.caravenna@math.unipd.it}
\and
Giambattista Giacomin \at Universit{\'e} Paris Diderot (Paris 7) and 
Laboratoire de Probabilit{\'e}s et Mod\`eles Al\'eatoires (CNRS U.M.R. 7599),
U.F.R. Math\'ematiques, Case 7012 (Site Chevaleret), 75205 Paris cedex 13, France, 
\email{giacomin@math.jussieu.fr}
\and
Fabio Lucio Toninelli \at Ecole Normale Sup\'erieure de Lyon,
Laboratoire de Physique and CNRS,
UMR 5672, 46 All\'ee d'Italie, 69364 Lyon Cedex 07, France,
\email{fabio-lucio.toninelli@ens-lyon.fr}}
%
%
\maketitle


\abstract{We review the literature on the localization transition for the class of polymers with random
potentials that goes under the name of {\sl copolymers near selective interfaces}. We 
outline the results, sketch some of the proofs and point out the open problems in the field. 
We also present in detail some alternative proofs that simplify what one can find in the literature.
}

\section{Copolymers and selective solvents}
\label{sec:0}

\subsection{A basic model}
In \cite{cf:GHLO}, T. Garel, D. A. Huse, S. Leibler and H. Orland introduced a simple model 
in order to look into how {\sl the statistical behavior of macromolecules can be strongly affected by  randomness in the physico-chemical properties of their constituents}. They aimed at 
a special class of macromolecules of linear type, the 
{\sl random hydrophilic-hydrophobic copolymers}, in a medium of water and oil,
separated by an interface. Such a polymer chain is just made up of monomers that differ for their affinity for water (or oil). The affinity is reduced to a real parameter that we call {\sl charge}:
the charge of the $j$-th monomer is denoted in \cite{cf:GHLO}
by $\zeta_j$ and, in mathematical terms,
$\{\zeta_j\}_{j=1,2,\ldots}$ is an IID sequence of Gaussian  random variables,
with given mean and variance. 
In order both to conform with the mathematical literature  and to generalize slightly
the problem we will write
$\zeta_j$ as $\go_j+h$,
where $h\in \R$ (or $h\ge 0$ as we will do next) and $\go=\{\go_j\}_{j=1,2,\ldots}$ is an  IID sequence of random variables (often referred to as the {\sl disorder}) such that
\begin{equation} \label{eq:M}
\M(t) \, :=\, \bbE\left[ \exp\left( t \go_1\right)\right]\, < \infty\, ,   
\end{equation}
for every $t \in \R$ and such that $\bbE \go_1=0$, 
$\bbE \go_1^2 =1$.
Apart for the larger class of charges that we allow and for notations, the Hamiltonian of the polymer model set forth in 
\cite{cf:GHLO}  is
\begin{equation} \label{eq:GHLO}
H_{N, \go,h}(S)\, :=\, 
 \sum_{n=1}^N \left( \go_n+h\right) \sign (S_n)
\, ,
\end{equation}
where $N$ is the length of the polymer and  $S=\{S_0, S_1,\ldots\}$ is a simple symmetric random walk trajectory ($S_0=0$, $\{S_{n+1}-S_n\}_{n=0,1,\ldots }$ is a sequence of independent identically distributed symmetric random variables that take only values $\pm 1$: the law
of $S$ is denoted by $\bP$ and we stress
that $\go$ and $S$ are independent). 

\begin{figure}[t]
\begin{center}
\medskip
\includegraphics[width = .95\columnwidth]{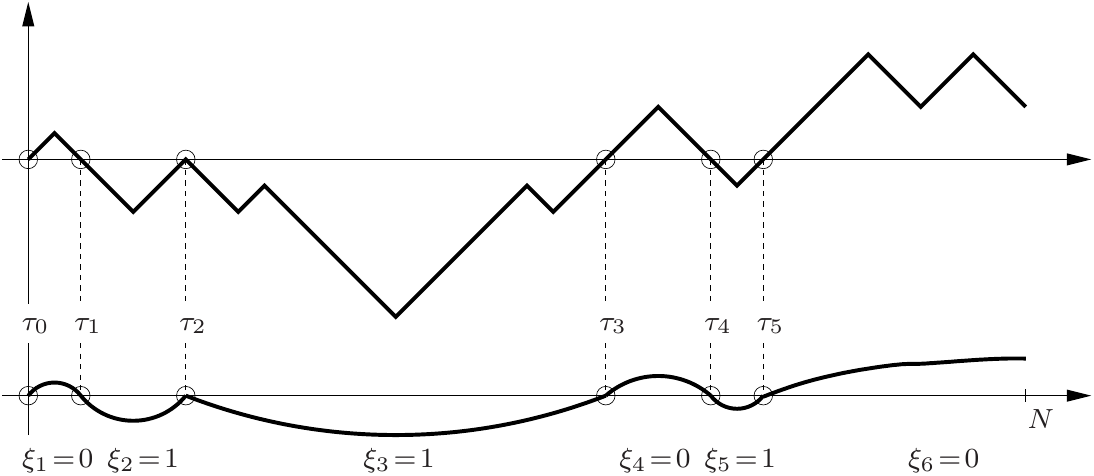}
\end{center}
\caption{\label{fig:basic}
The upper part of the figure shows a copolymer configuration: 
each {\sl bond}, or segment, $[(n-1, S_{n-1}), (n, S_n)]$ of the random walk trajectory
represents a monomer, so that $\sign(S_n)$ should be read as 
$+1$ (resp. $-1$) if this monomer is in the upper (resp. lower)
half-plane. Note that the Hamiltonian of the copolymer, cf. \eqref{eq:GHLO}, does not
depend on  the details of $S$ within an excursion, but only on 
the length and sign of the excursion. This naturally leads to 
the generalized model introduced in section~\ref{sec:generalized}
(lower part of the figure) in terms of a general 
discrete renewal process $\tau$ and a sequence of independent signs $\xi$. It is 
mathematically useful to use  the 
variable $\Delta_n$, that in the random walk case is just 
$(1-\sign(S_n))/2$, {\sl i.e.} the indicator function 
that the copolymer is below the axis.
Note in fact that, by \eqref{eq:GHLO},
$H_{N, \go, h} (S)= -2\sum_{n=1}^N
(\go_n+h) \gD_n + c_N(\go)$, with
$c_N(\go)=\sum_{n=1}^N (\go_n+h)$ which
does not depend on $S$ (or $\gD$), therefore
we can drop it without changing the polymer measure 
(that is what we do in \eqref{eq:taumodel}).}
\end{figure}

We invite the reader to have a look from now at Figure~\ref{fig:basic} for the directed polymer interpretation
of the trajectories of the model.
A small detail to deal with is $\sign(0)$: $\sign(S_n)$ should be read as
$\sign(S_{n-1})$ when $S_n=0$ and  this convention is particularly natural in directed polymer terms, 
because $\sign(S_n)$ is $+1$ ($-1$)
if the $n^{\textrm{th}}$ monomer is in the upper (lower) half plane, that is in oil (water),
see Figure~\ref{fig:basic}.
Still to conform with most of the mathematical literature on copolymers, the inverse temperature is denoted by
$\gl (\ge 0)$ instead of the more customary $\gb$, so that the Boltzmann factor that defines
the polymer model of length $N$ is $\exp(\gl H_{N, \go,h}(S))$. We are interested in the {\sl quenched} 
system so we underline the very different nature of the two sources of randomness: $\go$ is chosen 
once for all at the beginning of the experiment (the hydrophilic or hydrophobic character of the monomers 
does not change, while the chain fluctuates). 

At a superficial  level the effect of the charges on the polymer is quite intuitive: for $\gl>0$
positively charged monomers  ($\go_n+h>0$, that is hydrophobic monomers) prefer lying in the
upper half-plane (oil) and the opposite is true for the negatively charged ones.
But for large $N$ these {\sl energetically favorable} trajectories become more and more
atypical for $\bP$ since placing the monomers in their preferred solvent strongly
reduces the fluctuation freedom  of the chain. We are therefore dealing with an energy-entropy
competition that in the limit $N \to \infty$ leads to a localization-delocation transition: localization arises
when energy prevails and the polymer sticks to the oil-water interface,  visiting thus both oil and water,
while delocalization corresponds to the case in which
the polymer prefers to stay away from the interface. We will come back to this
with much more details, but we anticipate that this entropy-energy competition turns
out to be rather challenging: the phase diagram of the model is for the moment only
partially understood and  sound conjectures (or even only convincing physical heuristics) are lacking
on several fundamental issues. 
This may appear rather surprising in view of the very simple nature
of the model and of the fact, at the heart of the motivation of \cite{cf:GHLO}, 
that it represents one of the simplest instances of
a general mechanism that plays a crucial role in a variety of extremely important phenomena
(protein folding, to name one). 


\subsection{The (general) copolymer model}
\label{sec:generalized}


As argued in Figure~\ref{fig:basic} and its caption, the basic copolymer model does not depend on all
the details of $S$, but just on its zero level set, which is a renewal set, and on the sign of the excursions,
that is simply an independent fair coin tossing sequence. 
It is therefore natural, and at times really helpful, to look at the following generalized framework. Let us consider a
 discrete renewal process $\tau = \{\tau_n\}_{n\ge 0}$
on the non-negative integers $\N \cup \{0\}$, {\sl i.e.}, a sequence of random variables
such that $\tau_0=0$ and  $\{\tau_{j+1}-\tau_j\}_{j=0,1, \ldots}$ is an IID sequence of positive integer-valued variables, with marginal law satisfying
\begin{equation} 
\label{eq:K}
K(n)\, :=\, \bP \left( \tau_1 =n\right) \,
\overset{n\to\infty}{\sim}\, \frac{\C_K}{n^{1+\ga}}, 
\end{equation}
where $\ga >0$ and  $\C_K>0$
(we write $f_n \sim g_n$ for $f_n / g_n \to 1$).
Since we are dealing with renewal processes it is important to stress that $\bP(\tau_1< \infty )=1$,
that is $\sum_{n=1}^\infty K(n)=1$, so that 
 $\tau$ is a {\sl persistent} renewal. 
More generally, one could replace the constant $\C_k$ in \eqref{eq:K} by
a slowly varying function $L(n)$, but we stick for simplicity to the
purely polynomial asymptotic behavior \eqref{eq:K}.
It is well-known that the first return time to zero  of the simple symmetric random walk
$\inf\{n > 0:\, S_n = 0\}$
satisfies \eqref{eq:K} (restricted to the even integers,
due to the usual periodicity issue) with $\alpha = \frac 12$.
In particular, the basic model presented in the previous subsection
is a special case of the generalized copolymer model we are defining,
as it will be clear in a moment.

\begin{remark}
It is practical to switch freely from looking at $\tau$ as a sequence to
considering it a {\sl random set}, so for example $\vert \tau\cap [0, N]\vert$
is the number of renewals up to  $N$, or $n \in \tau$ is the event that
there exists $j $ such that $\tau_j=n$. For a comprehensive references on renewal processes see
 for example \cite{cf:Asm}.
 \end{remark}

The renewal $\tau$ identifies the polymer-interface contacts: we still need to know 
whether the excursion is above or below the axis.
For this let $\xi = \{\xi_n\}_{n\in\N}$ denote an IID sequence
of $B(1/2)$ variables (that is
$\bP(\xi_n = 0) = \bP(\xi_n = 1) = \frac 12$) independent of $\tau$, that we still
call signs.
Starting from the couple $(\tau,\xi)$  we 
build a new sequence $\gD = \{\gD_n\}_{n\in\N}$ by setting
$\gD_n = \sum_{j=1}^\infty \xi_j \, \ind_{(\tau_{j-1}, \tau_j]}(n)$,
in analogy with the simple random walk case:
 the signs $\Delta_n$ are constant between the points in $\tau$
and they are determined by $\xi$.

By 
 {\sl copolymer model} we mean 
the probability law $\bP_{N, \go} = \bP_{N, \go,\gl, h}$ for
the sequence $\Delta$ defined by
\begin{equation}
\label{eq:taumodel}
\frac{\dd \bP_{N, \go}}{\dd \bP} (\Delta) \, :=\, 
\frac 1{Z_{N, \go}}
\exp \left(-2 \gl \sum_{n=1}^N \gD_n (\go_n +h)
\right) \,,
\end{equation}
where $N \in \N$, $\gl\ge 0$, $h \in \R$ (but we can and will assume $h \ge 0$ without loss of generality) and $\go = \{\go_n\}_{n\in\N}$
has been introduced in \S~\ref{sec:0}.
The partition function
$Z_{N,\go} = Z_{N,\go,\gl, h}$ is given by
\begin{equation} \label{eq:discZ}
	Z_{N,\go} \;:=\; \bE \left[
	\exp \left(-2 \gl \sum_{n=1}^N \gD_n (\go_n +h)\right) \right] \,.
\end{equation}
In order to emphasize the value of $\ga$ in \eqref{eq:K},
we will sometimes speak of a {\sl $\ga$-copolymer model}, but
 $\bP_{N, \go}$ depends on the full distribution $K(\cdot)$, not only on $\ga$.

\subsection{The free energy: localization and delocalization}
\label{eq:fe}

We introduce the free energy of the 
copolymer by
\begin{equation} \label{eq:fe-qa}
	\tf (\gl, h) \,:= \, \lim_{N \to \infty } \, \tf_N(\gl, h) \,,
	\qquad \text{where} \qquad \tf_N(\gl, h) \,:=\,
	\, \frac 1N \, \bbE \left[ \log Z_{N, \go, \lambda ,h} \right] \,.
\end{equation}
The existence of such a limit follows by a standard argument
based on super-additivity, see for example
\cite{cf:dH} or
\cite[Ch.~4]{cf:Book}, where it is also proven that
\begin{equation}
\label{eq:fe-q}
\tf (\gl, h) \,= \, \lim_{N \to \infty }
\frac 1N  \log Z_{N, \go, \lambda ,h} \,, \qquad
\text{$\bbP(\dd \go)$-a.s. and in $L^1(\bbP)$}\,.
\end{equation} 
Equations \eqref{eq:fe-qa}--\eqref{eq:fe-q}
are telling us that the limit in  \eqref{eq:fe-q} does not depend
on the (typical) realization of $\go$, however it does depend
on $\bbP$, that is on the law of $\go_1$, as well as
on the inter-arrival law $K(\cdot)$.
This should be kept in mind, even if we omit
$\bbP$ and $K(\cdot)$ from the notation  $\tf (\gl, h)$.
Let us point out from now that $\tf(\gl, \cdot)$ and $\tf (\cdot, h)$ are convex functions,
since they are limits of convex functions. As a matter of fact $\tf(\cdot, \cdot)$ is only
separately convex because of the choice of the parametrization, but it is straightforward
to see that $(\gl, h) \mapsto\tf(\gl, h/\gl)$ is convex.

\begin{remark}
It is sometimes useful to consider the \emph{constrained}
partition function $Z_{N,\go}^\rc = Z_{N,\go,\lambda,h}^\rc$ of the model, defined by
\begin{equation} \label{eq:Zconstrained}
	Z_{N,\go,\lambda,h}^\rc := \bE \left[ \exp \left( -2\lambda \sum_{n=1}^N
	(\go_n + h) \Delta_n \right) \ind_{\{N \in \tau\}} \right] \,,
\end{equation}
which differs from \eqref{eq:discZ} only by the {\sl boundary condition}
factor $\ind_{\{N \in \tau\}} $.
It is a standard fact \cite[Remark~1.2]{cf:Book} that for all $N, \lambda, h$ we have
\begin{eqnarray} \label{eq:4}
	Z_{N,\go,\lambda,h}^\rc \le Z_{N,\go,\lambda,h} \le
	C\, N\, Z_{N,\go,\lambda,h}^\rc \,,
\end{eqnarray}
where $C$ is a positive constant. In particular, the free energy $\tf(\lambda, h)$
does not change if $Z_{N,\go,\lambda,h}$ is replaced by
$Z_{N,\go,\lambda,h}^\rc$ in \eqref{eq:fe-qa} and \eqref{eq:fe-q}.
Furthermore, since $N \mapsto \bbE(\log Z_{N,\go,\lambda,h}^\rc)$
is a real super-additive sequence, we can write
\begin{equation} \label{eq:feconstr}
	\tf(\lambda, h) = \lim_{N\to\infty} \frac 1N \, \bbE(\log Z_{N,\go,\lambda,h}^\rc)
	= \sup_{N\in\N} \frac 1N \, \bbE(\log Z_{N,\go,\lambda,h}^\rc) \,.
\end{equation}
\end{remark}

\smallskip

A crucial observation is:
\begin{equation} \label{eq:ebutc}
	\tf(\gl,h) \, \ge \, 0 \,  \text{ for every } \gl, h \ge 0 \,.
\end{equation}
This follows by restricting the expectation
in \eqref{eq:discZ} to the event $\{\tau_1 > N, \xi_1 = 0\}$,
on which we have $\Delta_1=0$, \ldots, $\Delta_N = 0$, hence we obtain
$Z_{N, \go} \ge \frac 12 \bP( \tau_1>N)$ and it suffices to
observe that $N^{-1}\log \bP( \tau_1>N) $ vanishes as $N \to \infty$,
thanks to \eqref{eq:K}. Notice that the event
$\{\tau_1 > N, \xi_1 = 0\}$ corresponds to the set of trajectories
that never visit the lower half plane, therefore the right hand side
of \eqref{eq:ebutc} may be viewed as the contribution to the free energy
given by these trajectories.

Based on this, we say that  $(\gl ,h) \in \cD$ ({\sl delocalized regime})
if $\tf(\gl,h)=0$, while $(\gl ,h) \in \cL$ ({\sl localized regime})
if $\tf(\gl,h)>0$. This may look at first as a cheap way to escape 
from the real localization/delocalization issue, that is inherently linked to
the path properties of the measure $\bP_{N, \go}$, but it is not the case.
Notice in fact that, if $h \mapsto \tf(\gl, h)$ is differentiable
(which fails at most for a countable number of values of $h$, by convexity),
by differentiating \eqref{eq:fe-q} and by convexity arguments we have 
\begin{equation}
\label{eq:path0}
-\frac 1{2\gl}
\frac{\partial}{\partial h} \tf(\gl,h)\, =\, \lim_{N \to \infty} \bE_{N, \go}\left[
\frac{\cN _N }N \right] \, , \qquad \bbP\text{-a.s.},
\end{equation}
where
\begin{equation}
\label{eq:cN}
\mathcal N_N\, :=\, \sum_{n=1}^N \Delta_n\, ,
\end{equation}
is just the total number of
monomers in the lower half-plane, that is in water
(cf. Figure~\ref{fig:basic}).
Therefore if $(\gl, h)$ is chosen in  the interior of $\cD$,
where $\tf \equiv 0$, the polymer visits water with null density
($\cN_N/N \to 0$). On the other hand, if $(\gl, h)$ is in $\cL$,
the polymer puts a positive density, precisely 
 $-(1/(2\gl))\partial \tf(\gl,h) / \partial h \in (0,1)$, of monomers in water and the rest, 
still a positive density, in oil.

We will deal below, {\sl cf.} section~\ref{sec:paths},
with sharper results on path behavior, but the elementary observation we have just made
shows that the definition we have set forth of localization and delocalization 
is far from being artificial. As a matter of fact,  it  is the natural physical definition,
and in fact it has been used already in \cite{cf:GHLO}, while in the mathematical literature 
was first introduced by \cite{cf:BdH}.


\subsection{The phase diagram}
\label{sec:phasediag}

Convexity and  the evident monotonicity of $\tf(\gl, \cdot)$ put
strong a priori constraints  on the phase diagram: let us go through 
this before going toward sharper questions. We can set 
\begin{equation}
h_c(\gl)\, :=\, \sup \{ h:\, \tf(\gl, h)>0 \}\, , 
\end{equation}
and the monotonicity of $\tf(\gl, \cdot)$ guarantees 
that $(\gl,h)\in \cL$ if and only if $h< h_c(\gl)$,
namely that
$h_c(\cdot)$ is the critical curve. 
Let us derive a number of elementary properties of
$h_c(\cdot)$.

The fact that 
$h_c(\gl) <\infty$ for every $\gl$ follows by the standard
{\sl annealed bound}:
\begin{multline}
\label{eq:annealing}
\bbE \log Z_{N, \go}\, \le \, \log \bbE Z_{N, \go}\, =\, \log \bE \bbE
\left[\exp \left(-2 \gl \sum_{n=1}^N \gD_n (\go_n +h)\right)\right]\, =\\
\log \bE \exp\left( (\log \M(-2\gl)-2 \gl h)  \sum_{n=1}^N \gD_n \right)\, , 
\end{multline}
so that $\tf(\gl, h)\le 0$ (hence, recall \eqref{eq:ebutc}, $\tf(\gl, h)= 0$) if  $h \ge \log \M(-2\gl)/(2 \gl)$,
and $\log \M(-2\gl)/(2 \gl) < \infty$ for every $\gl$ by \eqref{eq:M}. 

\begin{remark}
\label{rem:annealed}
The exponential  of the rightmost term in \eqref{eq:annealing}
is the partition function of the {\sl annealed} model associated to our quenched model. 
The free energy of the annealed model is rather trivial: it is in fact an elementary exercise to see that
\begin{equation}
\lim_{N \to \infty}\frac 1N \log \bE \exp\left( (\log \M(-2\gl)-2 \gl h)  \sum_{n=1}^N \gD_n \right)\,
=\, 
(\log \M(-2\gl)-2 \gl h)_+\, ,
\end{equation}
where $a_+:= a \ind_{a>0}$. 
The annealed model has therefore a (de)localization transition too:
its critical curve is $h_c^{ann}(\gl) := \log \M(-2\gl)/(2 \gl)$ 
and
we have just remarked that $h_c(\lambda) \le h_c^{ann}(\lambda)$
(this inequality on the critical curves is also referred to as annealed bound).
It can be noticed also that
the annealed free energy is not $C^1$ at criticality, that is the transition
is of first order. Let us stress that the annealed free energy looses a lot of details 
of the original model (in particular: no trace of  $K(\cdot)$!). 
\end{remark}

We have also $h_c(\gl)>0$ as soon as $\gl>0$, but 
this is not a trivial statement. In fact this means showing 
that $\tf(\gl, 0)>0$ for every $\gl > 0$: in more dramatic terms,
if the interaction does not select, on the average, a preferred solvent
($h=0$), the polymer is localized even at arbitrarily weak coupling, a result established first in \cite{cf:Sinai}. We skip the proof of this fact (in section~\ref{sec:loc} the proof of stronger results is  sketched) and we simply observe that, together with the annealed
upper bound, it implies that $h_c (\gl ) \longrightarrow 0$ as $\gl \searrow 0$.

At this point convexity  can be used in a very profitable way: since
$\{(\gl, y): \, \tf(\gl, y/ \gl) \le 0)\}$ is a convex set, its lower boundary
$\gl \mapsto \gl h_c(\gl)$ is 
a convex function. So we can write $h_c(\gl)=g(\gl)/ \gl$,
with $g(\cdot)$ convex such that $g(\lambda) = o(\lambda)$
as $\lambda \searrow 0$. This directly implies in particular continuity of $h_c(\cdot)$
and, with a little bit more of work, also the fact that 
$h_c(\cdot)$ is strictly increasing \cite{cf:BG}. 

\begin{figure}[t]
\begin{center}
\medskip
\includegraphics[width = .95\columnwidth]{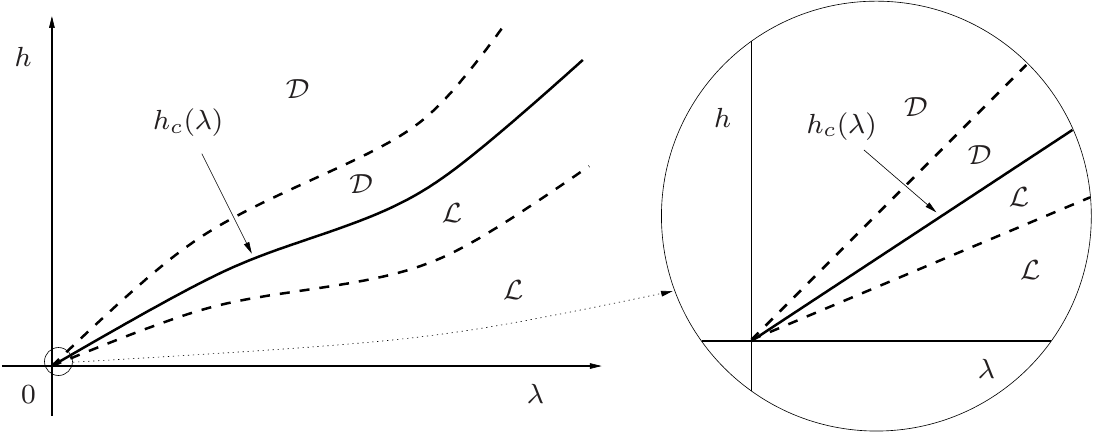}
\caption{\label{fig:phasediag} On the left there is a sketch of the phase diagram and the critical curve $\gl \mapsto h_c(\gl)$. 
The localized (resp. delocalized) regime corresponds to $(\gl,h)$ lying strictly below (resp. above) the critical 
curve. Explicit upper and lower bound on $h_c(\cdot)$ are known,
cf. Theorem~\ref{th:phasediag}, and they are schematically drawn as dashed lines. On the right one finds a zoom of the region near the origin, where the critical curve is close to 
a straight line and essentially all the
relevant information close to the origin (weak coupling regime) is encoded in the slope of this line: while the details of the critical curve do depend on the law of the disorder $\omega$ and on 
the details of $K(\cdot)$, the slope depends only on $\ga$. This is an important universality feature to which
 section~\ref{sec:BM} is devoted.  }
\end{center}
\end{figure}

Quite a bit of effort has been put into pinning down  the value of
$h_c(\cdot)$.
Figure~\ref{fig:phasediag} sums up the results that are known
on $h_c(\cdot)$ and, in particular, the content of 

\begin{theorem}
\label{th:phasediag}
For every $\lambda > 0$ the following explicit bounds hold:
\begin{equation}
\label{eq:LBUBhc}
\frac{1}{2\gl/(1+\ga)}\log \M\left( -2\gl/(1+\ga) \right)
\, \le \, h_c(\gl) \, < \, \frac{1}{2\gl}\log \M\left( -2\gl \right)
\,=\, h_c^{ann}(\lambda)\,,
\end{equation}
where the left inequality is strict when $\ga \ge 0.801$ (at least for $\gl$ small). 
\end{theorem} 

\noindent
The lower bound in \eqref{eq:LBUBhc} is proved in \cite{cf:BG}.
The strict inequality in the upper bound in \eqref{eq:LBUBhc} was first proved in \cite{cf:T_AAP}
to hold for large $\lambda$ and then extended to every $\lambda > 0$
in \cite{cf:BGLT}. In section~\ref{sec:deloc} we give an alternative, more direct proof.

Highlighted in Figure~\ref{fig:phasediag} is
the small $\gl$ behavior of $h_c(\gl)$. In fact, as more extensivily
explained in section~\ref{sec:BM}, for $\ga \in (0,1)$ we have 
$h_c(\gl) \stackrel{\gl \searrow 0}\sim m_\ga \gl$, with $m_\ga >0$ depending only on $\ga$. 
The slope $m_\ga$ is therefore a universal feature of the model: it does not depend
on the details of the disorder sequence $\omega$ and of the underlying renewal
$\tau$. The proof of such a result
goes through showing that for $\gl$ and $h$ small the free energy of the copolymer
is close to the free energy of a suitable continuum polymer model.

It becomes therefore quite relevant to get a hold of the value of 
$m_\ga$ (at least for $\ga \in (0,1)$). As a matter of fact,
from \eqref{eq:LBUBhc} one directly extracts 
\begin{equation}
\label{eq:dfg}
\frac 1{1+\ga} \, \le \, m_\ga \, \le \, 1\,, \qquad
\text{for every } \ga > 0 \, .
\end{equation}
but this result can be sharpened to
\begin{equation}
\label{eq:slopebounds}
\max\left(\frac 12, \frac{g(\ga)}{\sqrt{1+\ga}},\frac 1{1+\ga}\right) \, \le \, m_\ga \, < \, 1\, ,
\end{equation}
where $g(\cdot)$ is a  continuous function (of which we have
an  expression in terms of the primitive of an explicit function) such that
$g(\ga)=1$ for $\ga\ge 1$ and for which one can show 
that $g(\ga)/\sqrt{1+\ga} > 1/(1+\ga)$ for $\ga\ge 0.801$
(by evaluating $g(\cdot)$ numerically one can extend this result to $\ga \ge 0.65$). Since the existence of $m_\ga$ is not guaranteed for $\ga \ge 1$, to be precise both in
 \eqref{eq:dfg} and in  \eqref{eq:slopebounds}
 $m_\ga$ should be replaced by the inferior and superior limits of
 $h_c(\gl)/ \gl$. 
 
 \begin{remark} For sake of conciseness we have left aside the 
 $\ga =0$, that would require replacing the power law behavior \eqref{eq:K} with a regularly varying behavior
(allowing, in particular, the presence of logarithmic multiplicative corrections). The bounds
that we have just presented directly generalize \cite[Ch.~6]{cf:Book} and for $\ga=0$ the three terms in 
\eqref{eq:LBUBhc} coincide.
 \end{remark}

We would like to stress that Theorem~\ref{th:phasediag}
and the bounds \eqref{eq:slopebounds} show that
some claims in the physical literature are wrong. Notably 
in \cite{cf:GHLO,cf:MT} it is claimed, for $\ga=1/2$,  that $h_c(\cdot)$ coincides
with $h_c^{(ann)}(\cdot)$  (and \eqref{eq:slopebounds} shows that
also the weaker claim that $h_c (\gl) \stackrel{\gl \searrow 0}\sim
h_c^{(ann)}(\gl)$ is not correct. In \cite{cf:SSE,cf:Monthus} it is claimed 
that, still for $\ga=1/2$, the inequality in the left-hand side
of \eqref{eq:LBUBhc} is an equality. Theorem~\ref{th:phasediag}
falls short of proving that also this claim is false (even if it suggests it).
 A numerical study \cite{cf:CGG} lead for $\ga=1/2$, complemented by a careful 
statistical analysis using concentration inequalities, 
strongly suggests that the lower bound in  \eqref{eq:LBUBhc}
is strict and that the critical curve is somewhat {\sl halfway} between
the lower and the upper bound. 

What the previous results and discussion expose is 
that a convincing heuristic theory predicting the location of the critical 
curve, or just its slope at the origin, is lacking. In this sense 
we consider that capturing the value of $m_\ga$, for $\ga \in (0,1)$,
is an important open problem. Since computing explicitly quenched quantities may be really
out of reach, here are two sub-problems that are open:
\begin{itemize}
\item [$\star$] show that $m_{1/2}>2/3$;
\item [$\star$] is $m_\ga=1/(1+\ga)$ for some $\ga >0$?
\end{itemize}

Finally, we point out that a
{\sl reduced (simplified) copolymer model} was introduced in \cite{cf:BG},
taking inspiration from the approach in \cite{cf:Monthus}.
The original hope was that this simplified
model could catch the main features of the original copolymer model.
However, this does not seem to be the case, since it has been shown
\cite{cf:T_AAP,cf:BCT} that for the reduced model
one has $m_\alpha = 1/(1+\alpha)$ for all $\alpha \in (0,1)$.

\subsection{The critical behavior and a word about pinning models}
\label{sec:crit+pin}
Claims can be found in the physical literature about the critical
behavior of this model (at least in the original set-up, $\ga=1/2$, {\sl e.g.} \cite{cf:CW,cf:Habi,cf:Monthus,cf:MT}), but these 
claims do not always agree with each other, apart for the fact that they all
claim, not surprisingly, a smoothing effect of disorder. A rigorous result available on this issue
has been proven in 
\cite{cf:GT_cmp}: the 
transition of the general copolymer model is smooth (the derivative of the free energy vanishes at
least linearly when the critical point is approached, hence it is at least Lipschitz continuous 
at the critical point), in contrast with the annealed case (where the
the derivative of the free energy has a discontinuity at criticality, {\sl cf.} Remark~\ref{rem:annealed}):
for every $\lambda>0$ there exists $c(\lambda)<\infty$ such that
\begin{equation}
  \label{eq:5}
  \tf(\lambda,h_c(\lambda)+\delta)\le (1+\alpha)c(\lambda)\delta^2
\end{equation}
for every $\delta>0$. The result was obtained under some technical
conditions
on the disorder law, which are satisfied for instance in the case of
Gaussian or bounded charges. One can of course wonder whether 
the critical behavior of the copolymer model depends or not on $\ga$ (and on $\lambda$?)
and how,  but
once again, the substantial lack of sound 
physical predictions is quite disappointing. A natural open question however is:
\begin{itemize}
\item[$\star$]Can one improve \eqref{eq:5}, in the sense of replacing the exponent $2$ with a larger value?
\end{itemize}  
A somewhat deeper insight into this issue can be achieved by considering also another class of models, as we explain next.

The bound \eqref{eq:5} in fact coincides with the one available for disordered pinning 
models. Pinning models are a close companion to the copolymer,
since the Boltzmann factor of a pinning model is 
\begin{equation}
\exp\left( \sum_{n=1}^N (\gb \go _n + h) \gd_n\right)\, ,
\end{equation}
where $\gd_n = \ind_{n \in \tau}$, that is $\ind_{S_n=0}$ in the random walk set-up ($\go$ is chosen as before, so that $ \gb \go _n + h$
is a random variable of mean $h\in \R$ and variance $\gb^2$).
Therefore in this case the polymer has an interaction with the environment only when it touches the oil-water interface (or simply when
it touches the $x$ axis, usually called {\sl defect line}, 
since the model does not depend on the sign of the excursions). 
It is well known that pining models exhibit a localization transition
too and they can be dealt, to a certain extent, with similar techniques \cite{cf:Book}.
However, in the end, there are considerable differences, but let us try to single them out
and see what they suggest for copolymers:
\begin{enumerate}
\item The annealed pinning model is much richer than the annealed copolymer ({\sl cf.} Remark~\ref{rem:annealed}). In particular,
the annealed free energy does depend on $K(\cdot)$ and the critical behavior depends on $\ga$:
the transition is continuous (that is, the  free energy is $C^1$) as soon as $\ga \le 1$ and it becomes 
smoother and smoother as $\ga$ approaches $0$ \cite{cf:Book}. Harris criterion (see references in \cite{cf:Book,cf:cpam}) gives a precise prediction
 on what to expect for systems for which the annealed system
has a transition that is sufficiently smooth (for the pinning case the criterion boils down to $\ga <1/2$):
essentially it says that quenched and annealed systems have the same critical behavior and it
gives a precise prediction of the shift
in the critical point due to the disorder ({\sl irrelevant disorder regime}). At the same time it suggests/predicts  
that for $\ga>1/2$ {\sl disorder is relevant}, even if arbitrarily weak. This scenario has now been made rigorous,
see \cite{cf:cpam} and references therein, for pinning models. The crucial point for us is however the fact
that the free energy of the annealed copolymer model is not differentiable
at the critical point and therefore, in the Harris sense, disorder
is always relevant. 
\item Understanding critical phenomena when disorder is relevant is a major challenge and the possible scenarios 
set forth in the physical literature are quite intriguing, but very challenging and, at times, controversial (see {\sl e.g.} \cite{cf:MG,cf:Vojta} and references therein).
In this sense also the question that we have raised about improving \eqref{eq:5} acquires particular importance.
\item When $\ga>1$, also the annealed pinning model free energy is not $C^1$, and the critical curve
has been identified  with no more precision than for the copolymer model. In fact the annealed pinning critical curve
(again, a curve separating localized and delocalized regimes, in the $(\gb,h)$ plane)
behaves like
$ -\gb^2/2$ when $\gb $ is small and the quenched critical curve is in $[-c_+\gb^2,-c_-\gb^2]$ for $\gb $ small
(with explicit values of the constants $0< c_-<c_+ <1/2$,
cf. \cite{cf:DGLT}). This is absolutely parallel to the fact
that for the copolymer model $h_c(\gl) \in [c_- \gl, c_+\gl]$,
as one reads for example out of \eqref{eq:dfg},
\eqref{eq:slopebounds}. 
\end{enumerate} 

Finally, it is natural to wonder what happens when a pinning interaction is added to the
copolymer model, that is when not only the solvents are selective, but something special goes on
at the interface (thus taking into account  for example the lack of sharpness of the interface or the
fact that impurities could be trapped at the interface). There are works on this model, often called {\sl copolymer with adsorption} (see for example 
\cite{cf:Nicolas,cf:SW,cf:Whittington}), but the understanding is very limited: we refer to \cite[\S~6.3.2]{cf:Book}
for a detailed overview on this issue.


\subsection{Organisation of the paper}

The rest of the paper is devoted to going deeper into various results that we have stated, by 
giving either a sketch of arguments of proof, or alternative proofs and results that
complement what can be found in the literature. More precisely:
\begin{itemize}
\item In section~\ref{sec:loc} we give a sketch of the proof on the lower bound in~\eqref{eq:LBUBhc}.
\item In section~\ref{sec:deloc} we give an alternative proof of the upper bound in~\eqref{eq:LBUBhc}.
\item In section~\ref{sec:BM} we discuss the universality features of the
copolymer model in the weak coupling regime, i.e., for small values of $\lambda, h$.
\item Finally, section~\ref{sec:paths} is devoted to the description of the
available results on the path properties of the copolymer model.
\end{itemize}


\section{Localization estimates}
\label{sec:loc}

The aim of this section is to give a sketch of the proof of 
the lower bound in Theorem~\ref{th:phasediag},
that is of the left inequality in \eqref{eq:LBUBhc}, as well as of the left inequality in \eqref{eq:slopebounds}.

\medskip

The key-phrase for the approaches in this section is:
{\sl rare stretch strategies}. The idea, inspired by the renormalization group approach in 
\cite{cf:Monthus},
is to restrict the partition function to polymer trajectories that can visit the
lower half-plane only when there is a stretch of monomers that
are {\sl particularly, and anomalously,  hydrophilic}. To do this we introduce an intermediate scale $\ell$ (large, but fixed) and look
at the sequence of charges in blocks of $\ell$ charges at a time. 
We will consider two strategies (A and B) and, for simplicity, we will assume 
$\go_1\sim \cN(0,1)$:
\begin{itemize}
\item[A:] the $j^{\mathrm{th}}$ block is good if $\sum_{n=(j-1)\ell +1}^{j\ell} (\go_n+h) 
\le -m \ell$, with $m$ a positive value to be chosen below. For $\ell$ large, the probability that
a given block is good is very small, about $\exp(-\ell (h+m)^2/2)$, so that good blocks 
are typically separated by a distance of about   $\exp(\ell (h+m)^2/2)$.
\item[B:] the $j^{\mathrm{th}}$ block is good if $\sum_{n=(j-1)\ell +1}^{j\ell} (\go_n+h) 
= o(\ell)$. For $\ell$ large, the probability that
a given block is good is again very small for $h>0$, about $\exp(-\ell h^2/2)$, so that good blocks in this case 
are typically separated by a distance of about   $\exp(\ell h^2/2)$.
\end{itemize}

\begin{figure}[ht]
\begin{center}
\medskip
\includegraphics[width = .95\columnwidth]{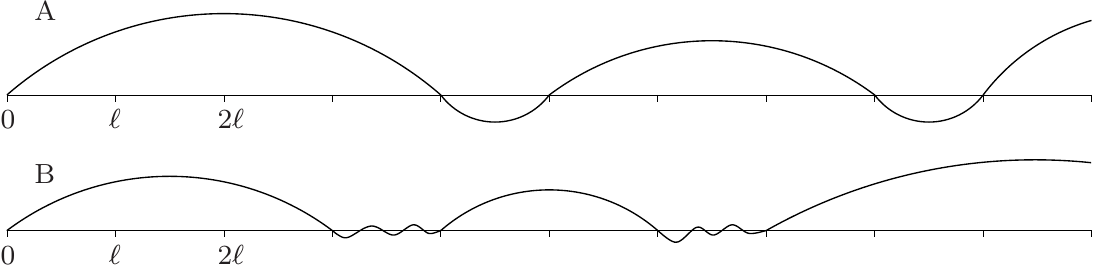}
\caption{\label{fig:strategies}
The lower bound strategies presented in this section are two and they are both
based on selecting some {\sl good blocks}, just by looking at the environment. They are actually blocks in which the environment is atypically negative: in strategy A we really select blocks in which the empirical average of the
charges $\go_n+h$ is smaller that a value 
$-m<0$ (selected in the end, in order to maximize the gain) and in strategy B we just aim at an empirical average close to zero (a rare event anyway, for $h>0$). Then for strategy
A we make a lower bound on the partition function by visiting the lower half plane if and only if a block is good
and by insisting that in a good block the walk stays in the lower half plane. In strategy B we target again 
the good blocks, but we put no constraint on the walk in the good blocks.}
\end{center}
\end{figure}

Given a sequence $\go$ of charges, the good blocks are identified
by the rules just given and we introduce
a set of polymer trajectories  $\gO_{N, \go}^A$, resp. $\gO_{N, \go}^B$,
for the strategy A, resp. B, defined as follows.  $\gO_{N, \go}^A$ is the set of trajectories 
that stay in the upper half plane except in presence of good blocks, when they stay 
in the lower half plane, see the upper part of Figure \ref{fig:strategies}.
The set of polymer trajectories  $\gO_{N, \go}^B$ still  includes 
all trajectories that stay in the upper half plane except in presence of good blocks, but when there is 
a good block the only restriction is that the polymer has to touch the oil-water
interface just before every good block and it has to touch it again at the end of the block,
see the lower part of Figure \ref{fig:strategies}.

The estimates are now just based on observing that 
\begin{multline}
\label{eq:d3f}
Z_{N, \go}\, \ge \, \bE \left[
	\exp \left(-2 \gl \sum_{n=1}^N \gD_n (\go_n +h)\right); \, \gO_{N, \go}^A  \right]
	\\
	=\, 
	\frac {K(i_1 \ell)}{2} \exp \left( -2 \gl \sum_{n=i_1\ell +1}^{(i_1+1)\ell }
	(\go_n +h) \right) \frac {K(\ell)}{2}\, \times \\ 
	\frac {K((i_2-i_1) \ell)}{2} \exp \left( -2 \gl \sum_{n=i_2\ell +1}^{(i_2+1)\ell }
	(\go_n +h)  \right) \frac {K(\ell)}{2} \, \times \, \ldots
\end{multline}
where the first good block is the $i_1^{\mathrm{th}}$ and so on. By using the definition of 
$\gO_{N, \go}^A$ we see that each term $-2 \gl \sum_{n=i_j\ell +1}^{(i_j+1)\ell }
	(\go_n +h) $ is bounded below by $2\gl m\ell$ and we notice that 
	the right-hand side of \eqref{eq:d3f} is a product of terms that, typically and to leading order,
	are just the same term, because the distance of good block is about
	$\exp(\ell (h+m)^2/2)$. Therefore this product of terms will
	give origin to an exponential growth in $N$ (localization!) if
	\begin{equation}
	\label{eq:g56}
-\frac{1+\ga}2 \ell (h+m)^2 + 2\gl \ell m +o(\ell)  	> 0
	\end{equation} 
(we have of course used \eqref{eq:K}). If we now optimize the choice of $m$ 
we readily see that this condition is met if $h< \gl/(1+\ga)$, for $\ell$ sufficiently large.
Therefore 	$h_c(\gl) \ge \gl/(1+\ga)$, which is the lower bound in \eqref{eq:LBUBhc} for
the special case of Gaussian charges (the extension to
general disorder is straightforward). 

Strategy B exploits a factorization similar to \eqref{eq:d3f} but this time the contribution 
in the good blocks is of about $\exp(\tf(\gl, 0) \ell)$, since in the good blocks 
the empirical average of the charges is zero, so the charges are essentially centered
and for large $\ell$ the random variables in a good block have a distribution that is close to 
centered IID standard Gaussian variables
(in the complete proof, {\sl cf.} \cite[Sec. 2]{cf:BGLT}, this step is performed via a change 
of measure argument that makes rigorous the heuristics presented here: the argument is slightly
more involved for non Gaussian charges).
The analog of \eqref{eq:g56} in this case becomes
\begin{equation}
	\label{eq:g56B}
-\frac{1+\ga}2 \ell h^2 + \tf(\gl, 0) \ell +o(\ell)  	\, . 
	\end{equation} 
At this point we need to estimate $\tf(\gl, 0)$. For example if 
one can show that $\tf(\gl, 0) \ge c \gl^2$ for come $c>0$ (say, for $\gl \le \gl_0$)
then for the same values of $\gl$ one would have 
that $h_c(\gl) \ge  \gl \sqrt{2c/(1+\ga)}$. This is the basic idea leading to
the (middle) lower bound in \eqref{eq:slopebounds}. Of course the work is now on estimating
$c$. We will not go into this issue which, ultimately, is a refinement of the result 
in \cite{cf:Sinai} and we refer to \cite[Sec. 2]{cf:BGLT} for details.
But we point out that
\begin{itemize}
\item as explained in section~\ref{sec:BM}, one can show that
$\tf(\gl,0)/ \gl^2$ has a positive limit that can be expressed in term of
the free energy of a continuum polymer;   
\item in order to improve on the lower bound in \eqref{eq:dfg} one needs
$\sqrt{2c/(1+\ga)}> 1/(1+\ga)$, that is $c>1/(2+2\ga)$. This can be established,
as recalled just below \eqref{eq:slopebounds}, for $\ga \ge 0.65$. 
If we were to improve on the lower bound in \eqref{eq:dfg} with this strategy 
for $\ga =1/2$ we would need to show $c>1/3$, but numerical estimations
suggest that $\lim_{\gl \searrow 0}\tf(\gl,0)/ \gl^2 $ is smaller (probably by little)
than $1/3$, so it is very likely that this strategy  (barely) fails to establish that 
the lower bound in \eqref{eq:dfg} can be made strict for $\ga=1/2$. The interest on this issue is because it would prove that the claims in \cite{cf:Monthus,cf:SSE} are not correct.
\end{itemize}


\section{Delocalization estimates}
\label{sec:deloc}

In this section we address the upper bound in Theorem~\ref{th:phasediag},
that is the right inequality in \eqref{eq:LBUBhc}. We will actually present (in full) an argument
that is substantially easier than the one originally used \cite{cf:BGLT}, even if it has the drawback to 
work only for $\ga \in (0,1)$. This argument is still based on the fractional moment method (first used in
the copolymer context in \cite{cf:T_AAP} to show that the upper bound in \eqref{eq:LBUBhc} holds for
$\gl$ large and unbounded disorder),
but it avoids the change of measure argument used in  \cite{cf:BGLT}. The change of measure argument in
\cite{cf:BGLT,cf:DGLT} is an important technique, as well as its refinement in \cite{cf:T_cg} that leads 
to the upper bound in \eqref{eq:slopebounds}, but we will not discuss these techniques here. 

\subsection{Fractional moment method: the general principle}
We consider the fractional moment method in its most elementary application.
For the constrained partition function \eqref{eq:Zconstrained} we can write
\begin{equation}
	Z_{N,\omega, \lambda, h}^\rc = \sum_{k=1}^N \sum_{0 = t_0 < t_1 < \ldots < t_k = N}
	\prod_{i=1}^k \frac{1 + e^{-2 \lambda
	(\sum_{t_{i-1}< j \le t_i} \omega_j + h(t_i -t_{i-1}))}}{2} \, K(t_i - t_{i-1}) \,.
\end{equation}
For $\gamma \in [0,1]$, from the inequality
$(a+b)^\gamma \le a^\gamma + b^\gamma$, valid for all
$a,b \ge 0$, we obtain the upper bound
\begin{equation} \label{eq:Zfrac}
	\bbE((Z_{N,\omega, \lambda, h}^\rc)^\gamma) 
	\le \sum_{k=1}^N \sum_{0 = t_0 < t_1 < \ldots < t_k = N}
	\prod_{i=1}^k \tilde K_{\gamma, \lambda, h}(t_i - t_{i-1}) \,,
\end{equation}
where we define
\begin{equation} \label{eq:tildeK}
	\tilde K_{\gamma, \lambda, h} (n) := \frac{1 + 
	e^{-(2 \lambda \gamma h - \log \M(-2\lambda \gamma)) n}}{2} \, K(n)^\gamma \,,
\end{equation}
We also set
\begin{equation}
	\Sigma(\gamma, \lambda, h) := \sum_{n\in\N}\tilde K_{\gamma, \lambda, h} (n) \,.
\end{equation}

Assume that $\Sigma(\gamma, \lambda, h) \le 1$.
By classical renewal theory \cite{cf:Asm}, the right
hand side of \eqref{eq:Zfrac} equals the probability that a renewal process
with step probability
(or sub-probability, if $\Sigma(\gamma, \lambda, h) < 1$)
$\tilde K_{\gamma, \lambda, h} (\cdot)$ passes through $N$;
in particular, it is bounded by $1$. Then
\begin{align*}
	\tf(\lambda, h) & = \lim_{N\to\infty} \frac{1}{N} \bbE(\log Z_{N,\omega, \lambda, h}^\rc)
	= \lim_{N\to\infty} \frac{1}{\gamma N} \bbE(\log (Z_{N,\omega, \lambda, h}^\rc)^\gamma) \\
	& \le \lim_{N\to\infty} \frac{1}{\gamma N} \log \bbE((Z_{N,\omega, \lambda, h}^\rc)^\gamma) = 0\,,
\end{align*}
whence $\tf(\lambda, h) = 0$ by \eqref{eq:ebutc}. 

This means that,
when there exists $\gamma \in [0,1]$ such that $\Sigma(\gamma, \lambda, h) \le 1$,
it follows that $(\lambda, h) \in \cD$.
This allows to give explicit estimates on the delicalized region.
Note that for $\gamma = 1$ we find the  {\sl  annealed
delocalized regime}, that we have already introduced: in fact
$\Sigma(1, \lambda, h) \le 1$ when
$h > h_c^{ann}(\lambda) := \log\M(-2\lambda)/(2\lambda)$.
Since for $\gamma \in [0, 1/(1+\alpha)]$ one sees immediately
that $\Sigma(\gamma, \lambda, h) = +\infty$,
the interesting range is $\gamma \in (1/(1+\alpha), 1)$.

\subsection{Fractional moment method: application}
Let us define for $\lambda > 0$
\begin{equation} \label{eq:hbar}
	\overline h(\lambda) := \inf \{h > 0: \ \exists \gamma \in [0,1]
	\text{ such that } \Sigma(\gamma, \lambda, h) < 1 \} \,.
\end{equation}

\begin{proposition}\label{th:frac}
For $\alpha \in (0,1)$
we have $h_c(\lambda) \le \overline h(\lambda) < h_c^{ann}(\lambda)$ for every $\lambda > 0$.
\end{proposition}

\noindent{\it Proof.}
We have just remarked that $\Sigma(\gamma, \lambda, h) \le 1$
implies $(\lambda, h) \in \cD$, therefore $h_c(\lambda) \le \overline h(\lambda)$.
It remains to show that $\overline h(\lambda) < h_c^{ann}(\lambda)$ for every $\lambda > 0$,
that is, for $\epsilon > 0$ sufficiently small
we can choose $\gamma \in [0,1]$
such that $\Sigma(\gamma, \lambda, h_c^{ann}(\lambda) - \epsilon) < 1$.
Note that for $\gamma \in (0,1)$
\begin{align*}
	& \frac{\partial \Sigma}{\partial \gamma} (\gamma, \lambda, h_c^{ann}(\lambda))
	= \sum_{n\in\N} \frac{1 + 	e^{-(\gamma \log \M(-2\lambda) - \log \M(-2\lambda \gamma)) n}}{2}
	(\log K(n)) \, K(n)^\gamma \\
	& \ - 2\lambda\left[ (\log \M)'(-2\lambda)
	+ \frac{\log \M(-2\lambda)}{2\lambda} \right] \sum_{n\in\N}
	\frac{e^{-(\gamma \log\M(-2\lambda) - \log \M(-2\lambda \gamma)) n}}{2} 
	\, n \, K(n)^\gamma \,.
\end{align*}
By the strict convexity of $\log\M(\cdot)$
and the fact that $\log\M(0)=0$,
\begin{equation}
	(\log \M)'(-2\lambda) < \frac{\log \M(0)- \log \M(-2\lambda)}{2\lambda}
	= - \frac{\log \M(-2\lambda)}{2\lambda} \,.
\end{equation}
Recalig our assumption \eqref{eq:K},
for $\alpha \in (0,1)$ we have $\sum_{n\in\N} n K(n) = \infty$,
therefore by Fatou's lemma
\begin{equation}
	\frac{\partial \Sigma}{\partial \gamma} (1^-, \lambda, h_c^{ann}(\lambda))
	:= \lim_{\gamma \uparrow 1} 
	\frac{\partial \Sigma}{\partial \gamma} (\gamma, \lambda, h_c^{ann}(\lambda))
	= + \infty \,.
\end{equation}
Since $\Sigma(1, \lambda, h_c^{ann}(\lambda)) = 1$, it follows that 
$\Sigma(1 - \eta, \lambda, h_c^{ann}(\lambda)) < 1$, for $\eta > 0$
small enough. By continuity, $\Sigma(1 - \eta, \lambda, h_c^{ann}(\lambda) - \epsilon) < 1$
for $\epsilon$ small enough, and the proof is completed.
\qed


\section{Continuum model and weak coupling limit}
\label{sec:BM}

In this section we explain in some
detail the universality feature sketched in Figure~\ref{fig:phasediag} and its caption.
The idea is that at weak coupling the details of the model, that is the law of the renewal beyond the exponent
$\ga$ and the law of the disorder, are inessential and a suitable continuum model 
captures the leading behavior of the (large class of) discrete models we consider.

As we already remarked, it is  convenient to look at the renewal process
$\tau = \{\tau_k\}_{k\ge 0}$ as a random subset of $[0,\infty)$.
It follows from our assumption \eqref{eq:K} that the rescaled random
set $\epsilon \tau = \{\epsilon \tau_k\}_{k\ge 0}$ converges in distribution
as $\epsilon \searrow 0$ toward a limit random
set $\tilde\tau^\alpha$, the so-called $\alpha$-stable regenerative set 
(we refer to \cite{cf:FFM} for more details; cf. also \cite{cf:CG2}).
This is a random closed subset of $[0,\infty)$
which is scale-invariant ($c \tilde\tau^\alpha$ has the same law as
$\tilde\tau^\alpha$, for every $c > 0$), has zero Lebesgue measure
and no isolated points.
For $\alpha = \frac 12$ we have the representation
$\tilde\tau^\alpha = \{t \in [0,\infty):  B_t = 0\}$,
where $B$ is Browian motion.

The complementary set $(\tilde\tau^\alpha)^c$, being open,
is the countable union of disjoint open intervals $\{I_n\}_{n\in\N}$.
We can then define a continuous-time process $\tilde\Delta^\alpha 
= \{\tilde\Delta^\alpha_t\}_{t\in [0,\infty)}$, which is constant
on each $I_n$ and takes the value $0$ or $1$, decided by
fair coin tossing: more precisely, in analogy with the discrete case, we set $\Delta_t :=
\sum_{n\in\N} \xi_n \, \ind_{I_n}(t)$ where
$\{\xi_n\}_{n\in\N}$ are i.i.d. $B(\frac 12)$ random
variables. For $\alpha = \frac 12$
we have the representation $\tilde\Delta^\alpha_t 
= \ind_{\{B_t < 0\}}$ with $B$ a Brownian motion. In general,
$\tilde\Delta^\alpha$ may be viewed as the limit in distribution of the
rescaled discrete process 
$\{\Delta_{\lfloor t/\epsilon \rfloor}\}_{t \in [0,\infty)}$
as $\epsilon \searrow 0$.

\smallskip

Let now $(\beta = \{\beta_t\}_{t \in [0,\infty)}, \bbP)$ be a Brownian motion,
independent of $(\tilde\Delta^\alpha, \bP)$.
We proceed for a moment in a somewhat informal way:
as $a \searrow 0$, $a^{-1} \omega_{\lfloor t/a^2 \rfloor}$
converges toward the white noise $\dd\beta_t/\dd t$ and
$\Delta_{\lfloor t/a^2 \rfloor}$ converges toward $\tilde\Delta_t$, therefore
\begin{gather*}
	a \lambda \sum_{n=1}^{N/a^2} (\omega_n + a h) \Delta_n
	= \lambda \int_0^N (a^{-1} \omega_{\lfloor t/a^2 \rfloor} + h)
	\Delta_{\lfloor t/a^2 \rfloor}\, \dd t
	\approx  \lambda \int_0^N (\dd \beta_t + h \dd t)
	\tilde \Delta^\alpha_t \,.
\end{gather*}
This hints at introducing a continuum partition function
\begin{equation} \label{eq:contZ}
	\tilde Z_{t,\beta} = \tilde Z_{t,\beta,\gl, h} \;:=\; \bE \left[
	\exp \left(-2 \gl \int_0^t \tilde\Delta^\alpha_s
	(\dd \beta_s + h \, \dd s)\right) \right] \, ,
\end{equation}
so that, recalling \eqref{eq:discZ}, one should have
$Z_{N/a^2, \omega, a \lambda, ah} \approx \tilde
Z_{N, \beta, \lambda, h}$ for $a$ small.

We now turn to precise statements. One can show that
the definition \eqref{eq:contZ} of the continuum partition function is well-posed,
for $\bbP$-a.e. $\beta$, and one introduces the corresponding continuum
free energy $\tilde\tf^\alpha(\lambda, h)$ in the usual way:
\begin{equation}\label{eq:contF}
	\tilde\tf_\alpha(\lambda, h) = \lim_{t \to \infty}
	\frac{1}{t} \, \bbE \log \tilde Z_{t,\beta,\gl, h} \,.
\end{equation}
The existence of this limit and the fact that it is self-averaging
(i.e., the expectation $\bbE$ can be dropped) require a much
longer and technical proof than the discrete counterpart, cf. \cite{cf:CG2}.

Also the continuum free energy
is non-negative: $\tilde\tf_\alpha(\lambda, h) \ge 0$
for all $\lambda, h \ge 0$, as one can easily check. The localized and delocalized regimes can
therefore be defined in analogy with the discrete case, namely
$(\lambda, h) \in \tilde\cL$ (resp. $(\lambda, h) \in \tilde\cD$) if
$\tilde\tf_\alpha(\lambda, h) > 0$ (risp. $\tilde\tf_\alpha(\lambda, h) = 0$),
and it is easily shown that they are separated by a critical
curve: $\tilde \cL = \{(\lambda, h) : h < \tilde h_c^\alpha(\lambda)\}$
and $\tilde \cD = \{(\lambda, h) : h \ge \tilde h_c^\alpha(\lambda)\}$.
There is however a major simplification
with respect to the discrete case: the scaling properties of
the processes $\tilde\Delta^\alpha$ and $\beta$ yield easily
$\tilde \tf_\alpha(a \lambda, a h) = a^2 \tilde \tf_\alpha(\lambda, h)$
for all $\lambda, h, a \ge 0$, therefore the critical curve is a straight
line: $\tilde h_c^\alpha(\lambda) =  m_\alpha \lambda$ 
for some $ m_\alpha$.

We can finally come back to
the rough consideration $Z_{N/a^2, \omega, a \lambda, ah} \approx \tilde
Z_{N, \beta, \lambda, h}$, that was discussed above.
This can be made precise in the
form of the following theorem.

\begin{theorem}
\label{th:main}
For an arbitrary discrete $\ga$-copolymer model we have
\begin{equation}
\label{eq:mainfe}
\lim_{a \searrow 0} \, \frac 1{a^2} \tf ( a\gl, a h)\, 
\;=\; \tilde \tf_\ga (\gl, h) \,, \qquad \forall \gl,h \ge 0 \,.
\end{equation}
Moreover 
\begin{equation}
\label{eq:mainslope}
\lim_{\gl \searrow 0} \frac{h_{c}(\gl)}{\gl}\, =\,  m_\ga.
\end{equation}
\end{theorem}

This result was first proved in \cite{cf:BdH} 
in the special case of the basic model of section~\ref{sec:0},
i.e., for the discrete copolymer model based on the simple random walk on $\Z$,
corresponding to $\alpha = \frac 12$
(in \cite{cf:GT} one can find an argument to relax the assumption in \cite{cf:BdH} of binary charges and in 
\cite{cf:Nicolas} the case with {\sl adsorption} is treated, {\sl cf.} the end of section~\ref{sec:crit+pin}).
The generalization to arbitrary
$\alpha$-copolymer models, with general disorder distribution,  is in \cite{cf:CG2}.

Note that \eqref{eq:mainfe} yields directly the existence of the
limit as $\lambda \searrow 0$ of $\tf(\lambda, 0)/\lambda^2$,
that was anticipated in section~\ref{sec:loc}, as well as the fact
that this limit coincides with $\tilde \tf_\alpha(1,0) > 0$.
We also point out that \eqref{eq:mainslope} is not a direct consequence of
\eqref{eq:mainfe}.

The importance of Theorem~\ref{th:main} relies in its \emph{universality} content:
for any fixed $\alpha \in (0,1)$ there is a single continuum model that captures
the behavior of all discrete $\alpha$-copolymer models for small values of $\lambda$
and $h$. In other words, the differences among these models become irrelevant
in the weak coupling limit. From this viewpoint, the slope $ m_\alpha$ of the
continuum critical curve is an extremely interesting object: improving the known
bounds $\frac{1}{1+\alpha} \le  m_\alpha < 1$ would yield a substantial
improvement in the understanding of the phase transition in this class of models.


\section{Path properties}
\label{sec:paths}

Up to now we have discussed the localization-delocalization
transition only in terms of free energy. A complementary, and equally
interesting, point of view is that of looking at path properties. In
other words, how does the typical (under $\bP_{N,\go}$, for typical
$\go$) polymer trajectory look like? The bottom-line of the picture
which has emerged up to this day is the following. In the localized
region the polymer makes order of $N$ excursions between the two
half-planes; the lengths of such excursions are $O(1)$ and their
distribution has an exponential tail.
In the delocalized region, on the other hand, the number of monomers in the defavorable
solvent (i.e. in the lower half plane) is not only sub-linear in $N$
(this information can be obtained immediately from the fact that the free
energy is zero there, cf. \eqref{eq:path0}) but it actually does not exceed
$O(\log N)$,
with high probability.
In the following, we discuss this picture in a bit more detail.

\subsection{The localized phase}

This subsection is is extracted from \cite{cf:GTloc}, to which we
refer for additional results, concerning for instance the exponential
tail of the length of the polymer excursions between the two solvents. Path properties 
in the localized phase have been studied also in \cite{cf:AZ,cf:BisdH}.

Let ${\bf M}_N:=\max_{i:\tau_i\le N}(\tau_i-\tau_{i-1})$ be the
length of the longest polymer excursion between the two solvents.
The following result says that  in the localized region
correlations decay exponentially fast, and the longest excursion is of
order $\log N$:
\begin{theorem}
\label{th:a1}
  Let $(\lambda,h)\in\mathcal L$. There exist constants $c_1,c_2$ such
  that, for every pair of bounded  local functions  $A,B$ of $\tau$ we have
  \begin{equation}
    \label{eq:8}
\sup_N    \bbE \left[|\bE_{N,\go}(A B)-\bE_{N,\go}(A)
      \bE_{N,\go}(B)|\right]\le c_1 \|A\|_\infty\,\|B\|_\infty e^{-c_2 d(A,B)},
  \end{equation}
where $d(A,B)$ denotes the distance between the supports of $A$ and
$B$.

Moreover, for every $\epsilon\in(0,1) $ the following holds in
$\bbP$-probability:
\begin{eqnarray}
  \label{eq:9}
  \lim_{N\to\infty}\bP_{N,\go}\left(\frac {1-\epsilon}{\mu(\lambda,h)}\le
  \frac{{\bf M}_N}{\log N}\le \frac {1+\epsilon}{\mu(\lambda,h)}\right)=1\,,
\end{eqnarray}
where $\mu$ is defined as
\begin{eqnarray}
  \label{eq:10}
  \mu(\lambda,h)=-\lim_{N\to\infty}\frac1N \log \bbE \frac{1+e^{-2\gl\sum_{n=1}^N(\go_n+h)}}{Z_{N,\go}}.
\end{eqnarray}
\end{theorem}
The existence of the limit \eqref{eq:10}, together with the bounds 
$0<\mu(\lambda,h)\le \tf(\lambda,h)$ in the localized phase ($\mu(\gl, h)\ge 0$ in general, and this is seen
in the same way as for $\tf(\gl, h)\ge 0$), is proven in \cite{cf:GTloc},
where one can also find an argument showing that $\mu(\gl, h) < \tf(\gl, h)$
under suitable (but not too restrictive) assumptions on the law of the charges. 

As a simple consequence of the exponential decay of correlations 
one can prove that 
\begin{enumerate}
\item the free energy is infinitely  differentiable (in
both $\lambda$ and $h$) in the localized region $\mathcal L$;

\item for every bounded local observable $A$ the limit
  $\lim_{N\to\infty} \bE_{N,\go}(A)   $ exists $\bbP( \dd\go)$
almost surely and is reached exponentially fast.
\end{enumerate}

As expected, the rate of exponential decay of the correlation functions (or inverse
correlation length), i.e. $c_2$
in \eqref{eq:8}, tends to zero as $h\nearrow
h_c(\lambda)$. In general, it is a very interesting open problem to
understand the relation between this and the way the free energy
vanishes
close to the critical point. In \cite{cf:T_jsp}, for the special case
where $K(\cdot) $ is the law of the
first return to zero of the symmetric simple
random walk on $\bbZ$, it was proven that the best constant $c_2$
coincides with $\mu(\lambda,h)$ defined in \eqref{eq:10}. 

There would be much to say about $\mu(\gl, h)$ and $\tf(\gl, h)$ and we prefer
to refer the reader to the introduction of 
\cite{cf:GTirrel} where this issue is treated in detail for pinning models. Here, in a simplistic way,
we just point out that 
Theorem~\ref{th:a1} and the bounds mentioned just after it are  telling us in particular 
that $\mu(\gl,h)$ is as good as  $\tf(\gl, h)$ for detecting the localization transition. But:
\begin{itemize}
\item [$\star$] is it true that $\log \mu(\gl,h) \sim \log \tf(\gl, h)$ as $h \searrow h_c(\gl)$?
In view of the discussion on disorder relevance in section~\ref{sec:crit+pin},
we expect that this is not the case and establishing such a result would be very interesting.
\end{itemize}

\subsection{The delocalized phase}

Recall the definition  $\mathcal N_N=\sum_{n=1}^N \Delta_n$ in \eqref{eq:cN}. The following theorem shows that,
strictly inside the delocalized region, $\mathcal N_N$ is typically at
most of order $\log N$:
\begin{theorem}[\cite{cf:GT}] \label{th:GT}
\label{th:2} For any $\delta>0,\lambda>0$ there exist $c>0,q>0$ such that for every
  $N\in \mathbb N$
  \begin{equation}
    \label{eq:1}
    \mathbb E\, {\bf P}_{N,\omega, \lambda,
    h_c(\lambda) + \delta}(\mathcal N_N\ge n)\le e^{-c\,n}
    \qquad \forall n\ge q\log N \,.
  \end{equation}
  
\end{theorem}

This result was proven in \cite{cf:GT} under the assumption that the
disorder law $\bbP$ satisfies a concentration inequality of
sub-Gaussian type. This holds for instance in the case of Gaussian or
bounded charges, and more generally whenever the distribution of
$\go_1$ satisfies a Log-Sobolev inequality.

Here we give a simpler argument, inspired by \cite{cf:Hubert}, which
works under the general assumptions of section \ref{sec:0} on the
disorder law and gives the weaker statement
\begin{eqnarray}
  \label{eq:2}
  \bbE \, \bE_{N,\omega, \lambda,
    h_c(\lambda) + \delta}(\mathcal N_N)\le \frac{c}{2\lambda\delta}\log N
\end{eqnarray}
for some constant $c$. The same argument also shows that at the critical
point $\cN_N$ is typically at most of order $\sqrt{N\log N}$:
\begin{eqnarray}
  \label{eq:3}
    \bbE \, \bE_{N,\go,\lambda,
    h_c(\lambda)}(\mathcal N_N)\le c'\sqrt{N\log N}
\end{eqnarray}
for some other constant $c'$.

Recalling \eqref{eq:discZ},
for all $h_0, h_1 \ge 0$ we can write
\begin{equation} \label{eq:basic}
	Z_{N,\go,\lambda,h_0}:=
	Z_{N,\go,\lambda,h_1}\cdot \bE_{N,\go,\lambda,h_1} \left[ \exp 
	\left( 2 \lambda (h_1 - h_0) \cN_N \right) \right] \,.
\end{equation}
Restricting $Z_{N,\go,\lambda,h_1}$ on the event $\{\tau_1 > N, \zeta_1 = 0\}$,
we obtain the bound $Z_{N,\go,\lambda,h_1} \ge \frac 12 \sum_{n > N} K(n)
\sim (const.) N^{-\alpha}$, by \eqref{eq:K}.
Applying Jensen's inequality we obtain
\begin{equation}
\begin{split}
	\frac 1N \bbE \left[ \log Z_{N,\go,\lambda,h_0} \right] & \ge
	\frac 1N \bbE \left[ \log Z_{N,\go,\lambda,h_1} \right] +
	\frac{2 \lambda (h_1-h_0)}{N} \bbE \, 
	\bE_{N,\go,\lambda,h_1}\left[ \cN_N \right] \\
	& \ge - c_1 \, \frac{\log N}{N} +
	\frac{2 \lambda (h_1-h_0)}{N} \bbE \, 
	\bE_{N,\go,\lambda,h_1} \left[ \cN_N \right] \,,
\end{split}
\end{equation}
where $c_1 \in (0,\infty)$. By \eqref{eq:4}
we have $\log Z_{N,\go,\lambda,h}^\rc \ge \log Z_{N,\go,\lambda,h}
- c_2 \log N$, for some constant $c_2$, therefore by \eqref{eq:feconstr}
we can write, for some $c \in (0,\infty)$,
\begin{equation} \label{eq:sa}
	\tf(\lambda, h_0) \ge \sup_{N\in\N} \left\{
	\frac{2 \lambda (h_1-h_0)}{N} \bbE \, 
	\bE_{N,\go,\lambda,h_1} \left[ \cN_N \right]
	- c \, \frac{\log N}{N}
	\right\} \,.
\end{equation}
\begin{itemize}
\item Take $h_0 = h_c(\lambda)$ and $h_1 = h_c(\lambda) + \delta$ with $\delta > 0$. Since
$\tf(\lambda, h_0) = 0$, from \eqref{eq:sa} we obtain \eqref{eq:2}.

\item Now take $h_0 = h_c(\lambda) - \delta$, with $\delta > 0$, and $h_1 = h_c(\lambda)$.
Since $\tf(\lambda, h_c(\lambda) - \delta) \le c_2 \delta^2$ for some
$c_2 = c_2(\lambda) \in (0,\infty)$ (cf. \eqref{eq:5}), it follows again from \eqref{eq:sa} that
\begin{equation}
	\limsup_{N\to\infty} \frac{\bbE \, \bE_{N,\go,\lambda,h_c(\lambda)}^\rc \left[ \cN_N \right]}
	{\sqrt{N \log N}} \le C  \,,
	\qquad \text{where} \quad C := \frac{\sqrt{c_1 \, c_2(\lambda)}}{\lambda} 
	\in (0,\infty) \,,
\end{equation}
whence \eqref{eq:3}.
\end{itemize}

\begin{remark}
Recalling \eqref{eq:tildeK}, let $\lambda, h > 0$ and $\gamma \in [0,1]$ be chosen such that
$\Sigma(\gamma, \lambda, h) < 1$, that is, $\tilde K_{\gamma, \lambda, h}(\cdot)$
is a sub-probability kernel on $\N_0$. Since $\tilde K_{\gamma, \lambda, h}(\cdot)
\sim c K(N)^\gamma \sim c' N^{-\gamma(1+\alpha)}$ as $N \to\infty$
for some constants $c, c' > 0$, by \eqref{eq:tildeK} and \eqref{eq:K},
it is a basic result in renewal theory that the right hand side of
\eqref{eq:Zfrac} is asymptotically equivalent to
$c'' \tilde K_{\gamma, \lambda, h}(N)$ as $N\to\infty$ for some
constant $c'' > 0$, cf. \cite[Theorem~A.4]{cf:Book}. Therefore
we have for all $N \in \N$
\begin{equation} \label{eq:Zgammac}
	\bbE((Z_{N,\omega, \lambda, h}^\rc)^\gamma) \le C_1
	\tilde K_{\gamma, \lambda, h}(N) \le C_1 \, K(N)^\gamma \,,
\end{equation}
where here and in the sequel $C_i$ denotes a generic positive constant.
Recalling \eqref{eq:4} and \eqref{eq:K}, 
for the original (non constrained) partition function we have
\begin{equation} \label{eq:Zgammaf}
	\bbE \left( (Z_{N,\omega, \lambda, h})^\gamma \right) \le C_2 \, N^\gamma \, K(N)^\gamma
	\le C_3 N^{-\gamma \alpha} \,.
\end{equation}
This relation can be exploited to improve Theorem~\ref{th:GT},
showing that with high probability $\cN_N$ is of order $1$. More precisely,
from the bound $Z_{N,\go,\lambda,h} \ge \frac 12 \sum_{n > N} K(n)
\sim (const.) N^{-\alpha}$ and \eqref{eq:basic} it follows that for any $\delta > 0$
\begin{equation}
	\bE_{N,\go,\lambda,h + \delta} \left[ \exp 
	\left( 2 \lambda \delta \cN_N \right) \right] 
	= \frac{Z_{N,\go,\lambda,h}}{Z_{N,\go,\lambda,h + \delta}}
	\le C_4 \, N^\alpha \, Z_{N,\go,\lambda,h} \,,
\end{equation}
hence applying \eqref{eq:Zgammaf} we obtain
\begin{equation}
	\bbE\left[ (\bE_{N,\go,\lambda,h + \delta} \left[ \exp 
	\left( 2 \lambda \delta \cN_N ) \right] \right)^\gamma \right]
	\le C_4^\gamma \, N^{\alpha \gamma} \, 
	\bbE \left( (Z_{N,\go,\lambda,h})^\gamma \right) \le C_5 \,.
\end{equation}
Recalling that $\gamma \in [0,1]$, by Markov's inequality we can write
\begin{align*}
	\mathbb E\, {\bf P}_{N,\omega, \lambda, h + \delta}(\mathcal N_N\ge n) & \le
	\mathbb E\, \left( {\bf P}_{N,\omega, \lambda, h+\delta}(\mathcal N_N\ge n)^\gamma \right) \\
	& \le e^{-2 \lambda \delta \gamma n} \, \bbE\left[ (\bE_{N,\go,\lambda,h + \delta} \left[ \exp 
	\left( 2 \lambda \delta \cN_N ) \right] \right)^\gamma \right]
	\le C_5 \,e^{-2 \lambda \delta \gamma n}\,. 
\end{align*}

Summarizing, we have shown that whenever relation \eqref{eq:Zgammaf} holds true,
there exists a constant $C > 0$ such that for every $\delta > 0$
\begin{equation}
	\mathbb E\, {\bf P}_{N,\omega, \lambda, h + \delta}(\mathcal N_N\ge n) \le
	C \,e^{-2 \lambda \delta \gamma n} \,, \qquad \forall n \in \N \,.
\end{equation}
We stress that relation \eqref{eq:Zgammaf} holds true
in particular for every $\lambda, h$ with
$h > \overline{h}(\lambda)$ (with a suitable choice of $\gamma \in [0,1]$,
recall \eqref{eq:hbar} and Proposition~\ref{th:frac}),
hence also \emph{below} the annealed critical curve. This is therefore an
improvement of equation (1.12) in~\cite{cf:GT}.
\end{remark}

We point out that delocalization properties were also studied in \cite{cf:BisdH}.
However the nature of the delocalized phase, in the pathwise sense, is 
still very little understood and, notably,  almost sure results are laking. For example:

\begin{itemize}
\item[$\star$] Is it true that, if $(\gl, h)$ is in the interior of $\cD$,  for every $\gep>0$ we have
$\lim_{N \to \infty} \bP_{N, \go} ( \gD_n=0$ for every $n \in [\gep N, N]\cap \N ) \, =\, 0$ $\bbP(\dd \go)$-almost surely?
\end{itemize} 



\begin{acknowledgement}
We gratefully acknowledge the support of the University of Padova
(F.C. under grant CPDA082105/08) and of ANR
(G.G. and F.L.T. under grant SHEPI).
\end{acknowledgement}



\end{document}